                  \def\version{14 October, 2020}		       %

\documentclass[reqno,12pt]{amsart} 
\usepackage[utf8]{inputenc}
\usepackage{amsmath} 
\usepackage[mathscr]{eucal}
\usepackage{amssymb}
\usepackage{srcltx} 
\usepackage{dsfont}
\usepackage[pagebackref, colorlinks=true,linkcolor=blue,citecolor=blue]{hyperref}
\usepackage{color}
\usepackage{tikz}
\usepackage{caption,subcaption}

\numberwithin{equation}{section}
 
\newfam\Bbbfam 
\font\tenBbb=msbm10 
\font\sevenBbb=msbm7 
\font\fiveBbb=msbm5 
\textfont\Bbbfam=\tenBbb 
\scriptfont\Bbbfam=\sevenBbb 
\scriptscriptfont\Bbbfam=\fiveBbb

\newcommand{\RR}     {\mathbb{R}} 
\newcommand{\ZZ}     {\mathbb{Z}} 
\newcommand{\NN}     {\mathbb{N}} 
\newcommand{\PP}   {\mathbb{P}} 
 
\newcommand{\EE}     {\mathbb{E}}

\def\1{{\mathchoice {1\mskip-4mu\mathrm l}      
{1\mskip-4mu\mathrm l} 
{1\mskip-4.5mu\mathrm l} {1\mskip-5mu\mathrm l}}} 
\newcommand{\ssup}[1] {{{\scriptscriptstyle{({#1}})}}} 
\def\comment#1{} 
\newtheoremstyle{thm}{2ex}{2ex}{\itshape\rmfamily}{} 
{\bfseries\rmfamily}{}{1.7ex}{} 
 
\newtheoremstyle{rem}{1.3ex}{1.3ex}{\rmfamily}{} 
{\itshape\rmfamily}{}{1.5ex}{}

\newtheorem{theorem}{Theorem}[section] 
\newtheorem{lemma}[theorem]{Lemma} 
\newtheorem{prop}[theorem] {Proposition}

\theoremstyle{definition}

\newcommand{\dd}{{{\rm d}}}

\newcommand{\cA}  {{\mathcal A}}

\newcommand{\cG}   {{\mathcal G }} 
 
\newcommand{\cI}   {{\mathcal I }}

\newcommand{\cM}   {{\mathcal M }} 
 
\newcommand{\cO}   {{\mathcal O }}

\newcommand{\cS}   {{\mathcal S }} 
\newcommand{\cT}   {{\mathcal T }}

\newcommand{\cX}   {{\mathcal X }}

\newcommand{\ee}   {{\operatorname e }}
\newcommand{\one}   {{\mathds{1}}}

\definecolor{Red}{rgb}{1,0,0}

\begin{document} 
 
\title[Branching random walks in random environment: a survey]{\large Branching random walks\\ \medskip in random environment: a survey}

\author[Wolfgang  K{\"o}nig]{} 

\maketitle

\thispagestyle{empty} 
\vspace{0.2cm}

\centerline {\sc By Wolfgang
K{\"o}nig\footnote{Technische Universit\"at Berlin, Str. des 17. Juni 136,
10623 Berlin, and Weierstrass Institute for Applied Analysis and Stochastics,
Mohrenstr. 39, 10117 Berlin, Germany, {\tt koenig@wias-berlin.de}. This text will appear as a chapter in the proceedings volume of the DFG Priority Programme 1590 {\em Probabilistic Structures in Evolution}}}
\medskip
\centerline{WIAS Berlin and TU Berlin}

\vspace{0.4cm}

\centerline{\small(\version)} 
\vspace{.5cm}

\begin{quote}
  We consider branching particle processes on discrete structures like
  the hypercube in a random fitness landscape (i.e., random
  branching/killing rates). The main question is about the location
  where the main part of the population sits at a late time, if the
  state space is large. For answering this, we take the expectation
  with respect to the migration (mutation) and the branching/killing
  (selection) mechanisms, for fixed rates. This is intimately
  connected with the parabolic Anderson model, the heat equation with
  random potential, a model that is of interest in mathematical
  physics because of the observed prominent effect of intermittency
  (local concentration of the mass of the solution in small
  islands). We present several advances in the investigation of this
  effect, also related to questions inspired from biology.
\end{quote}

\bigskip

{\em Keywords.} multitype branching random walk; random potential; parabolic Anderson model; Feynman-Kac-type formula; annealed moments; large deviations.

\medskip

{\em MSC2020.}  60J80; 60J55; 60F10; 60K37.

\section{Introduction}
\label{WK-sec:Intro}

In this chapter we study a topic that is of interest both in
mathematics and in mathematical population biology: branching random
walks\index{branching random walk}\index{branching!process}\index{branching!process!spatial} on a graph that either models the space in which we
live or a space of genotypes. That is,
the particles are under the influence of three random mechanisms:
movement, branching and killing. In the case of genotypes, the
movement can be understood as mutation\index{mutation!-selection!model},
and the killing as selection, hence such models belong to
the cornerstones of the mathematical description of random population
processes\index{population!process}. We introduce this mathematical
model in Section \ref{WK-sec:BRWRE} and explain its biological
interpretation in Section \ref{WK-sec:MutSel}.

The main point that we are interested in here is the situation where
both the branching rates and the killing rates may depend on the
state\index{branching!rates!state-dependent}\index{branching!rates!random}
that they are attached to, and that these rates are
taken as random, typically independent and identically distributed. In
this case, the expected number of particles in all the sites at a
given time is the solution to the heat equation with potential, which
is now random. This model is called the \textit{parabolic Anderson
  model} (PAM)\index{parabolic Anderson model} and has high physical
significance also in other parts of applied mathematics because of the
interesting effect of a high concentration property called
\textit{intermittency}\index{intermittency}. All this is explained in
Section \ref{WK-sec:BRWRE}.

In Section \ref{WK-sec:MultiType} we introduce a time-discrete version
of the model that has two state spaces: one describes the location,
the other describes the type of the particles. We develop a
version of a Feynman--Kac formula in this setting and derive the
large-time asymptotics for the expected total mass of the branching
process for a particular disctribution of the fitness variables. The
formulas that describe these asymptotics contain interesting
information about the counterplay of the two state spaces; mostly we
are interested in the description of the limiting structures that the
main part of the particles is sitting in.

A slightly different question is investigated in Section
\ref{WK-sec:hypercube}: For a particularly simple fitness distribution
on a finite but large graph, how much time is needed for the main part
of the population to move to the fittest site?

In Section \ref{WK-sec:HigMom}, we extend the setting of the parabolic
Anderson model (which describes the first moment of the number of
particles) to a Feynman--Kac type formula for all the moments of the
number of particles, taken over migration, branching and killing,
valid for all choices of the branching/killing rates. This formula is
employed afterwards to derive a large-time asymptotics for the high
moments, for a particular potential distribution.

In Section \ref{WK-sec:future} we mention and motivate future work on
extensions of the study of branching random walks in random
environment, like the study of the high-moment asymptotics, the
parabolic Anderson model on random tree-like graphs and the effect
coming from inserting a mutually repellent force between the
particles.

\section{Spatial branching random walks in random
  environment}\label{WK-sec:BRWRE}

We introduce the basic model in Section \ref{WK-sec:model},
explain its interpretation for evolutionary applications in Section
\ref{WK-sec:MutSel}, its relation with the parabolic Anderson model in
Section \ref{WK-sec:PAM} and give a brief account on the types of
limiting statements available in the literature in
Section \ref{WK-sec:MomAsy}.

\subsection{The model}
\label{WK-sec:model}

Let $\cX$ be some graph and $X=(X_t)_{t\in[0,\infty)}$ be a
Markov process in continuous time on $\cX$, which is most often taken
as the simple random walk. The main examples for the state space $\cX$
are the infinite space $\ZZ^d$ or the finite hypercube
$\cX =\{1,\dots,K\}^n$ for some $K,n\in\NN$. Let
$\xi_+=(\xi_+(x))_{x\in\cX}$ and $\xi_-=(\xi_-(x))_{x\in\cX}$
be two collections of nonnegative numbers, attached to the sites of
the state space. Now the spatial branching random walk\index{branching random walk!spatial} with potential $(\xi_+,\xi_-)$ is introduced as
follows. Starting with a single particle at time zero at the origin
$\cO\in\cX$, the particle runs like a copy of $X$ through
$\cX$. Furthermore, when located at $x$, the particle splits into two
particles with rate $\xi_+(x)$ and dies (is removed from the system)
with rate $\xi_-(x)$. Hence $\xi_+(x)$ is the \emph{branching rate} at
$x$, and $\xi_-(x)$ is the \emph{killing rate}. The two newborn copies
proceed independently in the same way as the initial particle.

This is one of the basic models for the time-evolution of a
population in a space with several random mechanisms: \emph{migration},
\emph{branching} and \emph{killing}. The fields of rates, $\xi_+$ and
$\xi_-$, introduce disorder in the latter two mechanisms; indeed we
will consider them here often as being random and independent and
identically distributed. Areas with high values of $\xi_+$ will
produce many offspring over large time, and areas with high values of
$\xi_-$ will contain fewer particles over long time. Since the $\xi_+(x)$
are rates, the number of particles at $x$ at time $t$ is likely to
grow exponentially fast with rate $\xi_+(x)$ (neglecting the killing
mechanism), but the migration mechanism diminishes this effect, since
the many newborn particles have the tendency to move away from this
site (but may return later). See \cite{WK-CM94} for more explanations
about random walks with branching and killing in a random medium, in
particular in connection with equations of the type that we will be
examining in Sections \ref{WK-sec:MutSel} and \ref{WK-sec:PAM}.

The space $\cX$ is often taken as (a discrete version of) the real
space in which we live, i.e., the Euclidean space $\ZZ^3$, and many
investigations of the model concentrate on this
interpretation. However, from the viewpoint of evolution of a
population in a biological sense, it is also highly interesting to
view $\cX$ as a space of genotypes, i.e., as a space of more abstract,
biological properties of the particles. One particularly descriptive
example is the choice of $\cX$ as the set of \emph{gene sequences},
representing each particle by its genome, which is a sequence composed
of four or, for simplicity, two alleles. Here the hypercube
$\cX =\{-1,1\}^N$ is the most canonical choice. In Section
\ref{WK-sec:MutSel} we explain the interpretation of our model from
the viewpoint of biological evolution.

The standard reflex of a probabilist is to choose the migration
mechanism as the nearest-neighbour simple random walk, where two gene
sequences are considered neighbours if they differ in precisely one
allele. We will consider this in Section \ref{WK-sec:BRWgraphs}, but
let us note that there are more biologically inspired migration
mechanisms that respect the way in which mutation and selection of
gene sequences really happen~\cite{WK-BG07}. A mathematical
investigation of the above branching model with this type of random
motion is widely open and would be rather interesting and
well-motivated.

A spatial branching process model (called the {\em logistic branching
random walk}) without randomness in the birth rates, but a kind of birth control in
the case of a local overpopulation is considered in \cite[Section~3]{WK-MBNG20} in this volume.
Here the birth rate in a given site
is bounded, and is small (possibly zero) if the current number of
individuals in a neighbourhood
is large. The main results reported there are criteria for survival
and for long-time distributional
convergence to equilibrium, conditioned on survival. Proofs are based
on a comparison
technique using coupling.

\subsection{The heat equation and a mutation-selection population
  model}\label{WK-sec:MutSel}

The model of branching random walks introduced in Section
\ref{WK-sec:BRWRE} on the hypercube $\cX=\{-1,1\}^N$ with random
branching rates has a profound interpretation in terms of a
mutation-selection population model in a random fitness
landscape\index{fitness!landscape!random}. (See \cite{WK-K19} for theoretical
background about the notion of fitness in biological context and about
the numbers and lengths of paths between sites of different
fitnesses.) In \cite{WK-AGH16}, this is described as
follows. \index{parabolic Anderson model!biological interpretation}

The mutation-selection model\index{mutation!-selection!model} is given
by the solution $v_N(t,\cdot,y)$ of the partial differential equation
(PDE)
\begin{equation}\label{WK-mutaselec}
  v^{}_N(t,x,y) =\frac{1}{N} \Delta v^{}_N(t,x,y) +[\xi(x)-\overline\xi(t)]\,
  v^{}_N(t,x,y),  \quad\ \ t\in[0,\infty), x\in \{-1,1\}^N,
\end{equation}
with the localised initial condition
$v_N(t,\cdot,y) =\delta_y(\cdot)$. The potential $\xi=\xi_+-\xi_-$ is
interpreted as a \emph{fitness landscape}, and the mean
fitness\index{fitness!mean} is given by
\begin{equation}\label{WK-meanfitness}
 \overline\xi(t)=\!\!\sum_{x\in\{-1,1\}^N}v_N(t,x,y) \xi(x).
\end{equation}
Here $\Delta$ is the \emph{Laplace operator}\index{Laplace!operator}
defined by
\begin{equation}\label{WK-LaplaceOp}
\Delta f(x) =\!\! \sum_{y\in \cX\colon y \sim x} \big( f(y) - f(x)\big)\, .
\end{equation}
Now, \eqref{WK-mutaselec} is not a particle model, but a PDE; however,
let us already note that such types of PDEs are very suitable for
describing the expectation of the number of particles in a spatial
branching process for fixed branching rates. This will be explained in
a broader perspective in Section \ref{WK-PAM}.

Let us briefly explain the biological meaning of
\eqref{WK-mutaselec}. Haploid genotypes are identified with linear
arrangement of $N$ sites $x= (x(1),\dots,x(N))$ with each site taking
values $-1$ or $+1$. In the multilocus context, sites correspond to
loci and the variables $x(i)$ to alleles. In the context of molecular
evolution, $x$ corresponds to a DNA (or RNA) sequence, where the
nucleotides are lumped into purines (say, $+1$) and pyrimidines
(say,$-1$). In the biology literature, the hypercube\index{hypercube!as sequence space} is usually called the sequence space. Then the
mutation-selection model in \eqref{WK-mutaselec} describes the
evolution of an infinite population of haploids that experience only
mutation and selection. The population evolves in continuous time
(non-overlapping generations) with mutation and selection occurring
independently (in parallel). $\xi(x)$ is the Malthusian fitness of type
$x$ and the $\xi (x)$ form a fitness landscape, which in our case is
random. Site mutations happen with rate $1/N$ (hence, at a total rate
of one). From \eqref{WK-meanfitness} it follows that
$\sum_{x} v_N(t,x,y) = 1$, and
$v_N(t,x,y)$ corresponds to the frequency of type $x$ under this
evolution. Finally, note that the localised initial condition means
that initially the population consists of only type $y$.

The competition between diffusion and potential in the PAM translates
into a competition between mutation and selection, two driving forces
of Darwinian evolution.  We refer to the classical book \cite{WK-CK70}
for an introduction to population genetics and to \cite{WK-GB99} for
an excellent survey that involves the statistical physics methods used
to solve mutation-selection models for a wide range of landscapes.

In the following, we will concentrate on the situation where the
fitness landscape\index{fitness!landscape!random} is random. The
motivation for this is the following. Realistic landscapes are
expected to be complex with structures such as valleys and hills.
Random fitness landscapes naturally form a class of complex
landscapes. The first obvious choice, that is, an i.i.d. landscape, is
also known as the \emph{house of cards} model\index{house of cards model}.

Let us also mention that it is well-known \cite{WK-M76} that
\begin{equation}
v^{}_N(t,x,y)=\frac{u^{}_N(t,x,y)}{u^{}_N(t,x)},
\end{equation}
where $u_N$ is the solution to \eqref{WK-mutaselec} with
$\overline\xi$ replaced by zero, which is called the \emph{heat
  equation} that we are going to introduce in Section
\ref{WK-sec:PAM}. In a way, $u_N(t,x,y)$ can be thought of as an
absolute frequency.

\subsection{The parabolic Anderson model}\label{WK-sec:PAM}

Let us return to the branching random walk model introduced
in Section \ref{WK-sec:model} with fixed branching rates $\xi_+(x)$
and killing rates $\xi_-(x)$. We want to explain this model from the
viewpoint of PDEs, like in \eqref{WK-mutaselec}. For the time being,
it is inessential that the rates are random. By $\eta(t,x)$ we denote 
the number of particles at time $t$ at site $x$, and by $u(t,x)$ 
its expectation with respect to migration, branching and
killing. Note that the expectation of $\eta(t,x)$ does not depend on
$\xi_+$ and $\xi_-$, but only on the difference field
$\xi=\xi_+-\xi_-$, since the branching and killing balances out in the
expectation. We can therefore view $u(t,\cdot)$ as a time-dependent
field that depends on $\xi=(\xi(x))_{x\in\cX}$.

As a matter of fact, the time-dependent field $u(t,\cdot)$ satisfies
the \emph{heat equation with potential} $\xi=\xi_+-\xi_-$,
\index{heat equation with potential} which is the equation
\begin{equation}\label{WK-PAM}
\begin{aligned}
\partial^{}_t u(t,x) &=  \Delta u(t,x) + \xi(x) u(t,x),
\qquad x  \in \cX , \, t \in (0,\infty),
\\
u(0,x) &=\one_0(x),\qquad\qquad\qquad\qquad\quad\! x \in \cX .
\end{aligned}
\end{equation}
Here $\Delta$ is the generator of the Markov process $X$, i.e., in
case of a simple random walk on the graph $\cX$, it is the
Laplace operator defined in \eqref{WK-LaplaceOp}.

The solution theory for the heat equation is rather rich and
explicit. Indeed, the solution $u$ to \eqref{WK-PAM} can be
represented in terms of the \emph{Feynman--Kac
  formula}\index{Feynman--Kac!formula}
\begin{equation}\label{WK-FK}
u(t,x)=\EE^{}_{\cO}\Big[ \ee^{\int_0^t \xi(X_s)\,\dd s}\one\{X^{}_t=x\}\Big],
\end{equation}
and in terms of an eigenvalue expansion\index{eigenvalue expansion}
for the eigenvalues $\lambda_1>\lambda_2\geq\lambda_3\geq\dots$ of the
operator $\Delta+\xi$,
\begin{equation}\label{WK-eigenexp}
u(t,x)=\sum_{k}\ee^{t\lambda^{}_k}\phi^{}_k(\cO )\phi^{}_k(x),
\end{equation}
where $(\phi_k)_k$ is an orthonormal sequence of corresponding eigenfunctions.

So far, the field $\xi$ of differences of branching and killing rates
collects just the coefficients of the PDE in \eqref{WK-PAM}, and their
randomness is inessential for the above. But let us turn now to the
case that $\xi$ is a \emph{random} collection of numbers, which we then
call a \emph{random potential}. In this case, the operator $\Delta+\xi$
and its spectrum are also random; it is called the \emph{Anderson
  operator}\index{Anderson!operator} and plays an important role in
mathematical physics; in particular the famous phenomenon of {\rm
  Anderson localisation}\index{Anderson!localisation} has kept the
interest in the description of its spectrum awake for decades. The
heat equation with random potential in \eqref{WK-PAM} is then called
the \emph{parabolic Anderson model}\index{parabolic Anderson model}.

As we are interested in the behaviour of $u(t,\cdot)$ for large $t$,
we will be concerned only with the boundary part of the spectrum,
mostly only with the principal eigenvalue $\lambda_1$ (this is easily
seen from \eqref{WK-eigenexp}). Anderson localisation is, roughly
speaking, the phenomenon that all the eigenfunctions of the Anderson
operator $\Delta+\xi$, at least the ones corresponding to eigenvalues
close to the boundary of the spectrum, are highly concentrated in
small areas that are randomly located. The principal eigenvalue is
included in this statement. Again looking at \eqref{WK-eigenexp}, we
can guess that this concentration property is inherited by the random
function $u(t,\cdot)$. This has been shown to be true in some
generality, and this is the basis of a detailed understanding of the
large-time behaviour of the entire branching process in a random
potential of branching/killing rates.

The main assumption on the random potential $\xi$ is that it is a
collection of independent, identically distributed (i.i.d.) random
variables. Many distributions can be considered, and for many
applications there is no canonical choice. Not even ubiquitous
distributions like the standard normal distribution or the exponential
distribution can be considered canonical from the viewpoint of
applications nor from the viewpoint of interesting emerging effects in
this model. Rather, it turned out in the study of the moments of $u$
that the most natural choice in this respect is perhaps the
\emph{double-exponential distribution}\index{distribution!double-exponential},
\begin{equation}\label{WK-DE}
  {\rm Prob} \big( \xi(0) > u \big) = \ee^{-{\rm e}^{u/\varrho}}, \qquad
  u \in \RR ,
\end{equation}
since it turned out in several investigations that the emerging
picture in the long-time behaviour of the process has a non-trivial
and not too complicated structure.

Our main interest here is in the description of the long-time (i.e.,
large-$t$) behaviour of the random field $u(t,\cdot)$, in particular
of its \emph{total mass}\index{branching!process!total mass}
\begin{equation}\label{WK-totalmass}
U(t)=\sum_{x\in\cX} u(t,x).
\end{equation}
Of particular interest is the question whether $u(t,\cdot)$ develops
some (randomly located) areas with particularly high values, i.e.,
particularly high accumulation of particles. Such areas are often
called \emph{intermittent islands}\index{intermittency}. Further
questions concern the size and the location of these islands and
characterisations of the potential $\xi(\cdot)$ and the solution
$u(t,\cdot)$ inside them.

\subsection{Asymptotics}\label{WK-sec:MomAsy}

Let us explain some of the fundamental results on the
large-$t$ behaviour of the parabolic Anderson model introduced in
Section \ref{WK-PAM}. See \cite{WK-K16} for a comprehensive survey on
the research until 2016 on the PAM, which was almost exclusively done
on the Euclidean space $\RR^d$ or $\ZZ^d$. We abstain in
this section from giving explicit references to the original
literature and refer to \cite{WK-K16} across-the-board.

Here is a fundamental assertion that has been proved about the
large-$t$ behaviour of the PAM: the asymptotics of the moments of the
total mass of the solution, $\langle U(t)\rangle$ (we write
$\langle\cdot\rangle$ for expectation with respect to the potential
$\xi$). We present this here only for the double-exponential
distribution, as this is the one that we will consider further. The
following is also the historically first result of this type; together
with its proof, it gave guidance to the analysis of the PAM for many
other potential distributions over the last 20 years. Introducing
$H(t)=\log\langle \ee^{t \xi(0)}\rangle$, the cumulant generating
function of one of the potential variables, we have the following (see
Theorem 3.13 and Remark 3.17 in \cite{WK-K16}).

\begin{theorem}[Moment asymptotics of the PAM]\label{WK-ThmMomAsy}
  We consider the PAM as in \eqref{WK-PAM} on $\cX=\ZZ^d$ and assume
  that the random potential $\xi=(\xi(x))_{x\in\ZZ^d}$ is i.i.d.\ and
  doubly-exponentially distributed as in \eqref{WK-DE} with parameter
  $\rho\in(0,\infty)$. Then there is a number
  $\chi=\chi(\rho)\in[0,2d]$ such that
\begin{equation}\label{WK-MomAsy}
  \frac{1}{t} \log \langle U(t)\rangle =\frac{H(t)}{t}-\chi+o(1),\qquad
  t\to\infty.
\end{equation}
\end{theorem}

Note that $H(t)=\rho t\log t+\rho t+ o(t)$ for $t\to\infty$, hence the
leading term $H(t)/t$ tends to infinity, while the second term $\chi$
is a constant.

Indeed, the characteristic quantity $\chi$ is represented by a certain
variational formula that is interpretable and can be investigated more deeply.
The interpretation is that the best contribution to the
expected total mass comes from a localised peak in a certain island in
$\ZZ^d$, called the \emph{intermittent island}\index{intermittency}.
More precisely, it comes from a pretty
small area in which the potential $\xi$ assumes particularly high
values. This area is essentially a centered ball with a deterministic,
fixed radius (not depending on $t$), which is a special feature that
only the double-exponential distribution shows. We can also argue in
terms of the Feynman--Kac formula in \eqref{WK-FK} that the random walk
$(X_s)_{s\in [0,t]}$ spends practically all the time in this island in
order to benefit as much as possible from the high potential
values. In terms of the eigenvalue expansion \eqref{WK-eigenexp}, we
can argue that the local principal eigenvalue (say, with zero boundary
condition) in that island is enormously large. The leading term
$H(t)/t$ expresses the high potential value, and the second term
$\chi$ expresses the optimal compromise between the long stay of the
random walk in the island and the size of the island. (In order to see
that the large-$t$ asymptotics may have anything to do with any kind
of optimisation, note that both the Feynman--Kac formula and the
eigenvalue expansion reveal that $U(t)$ is something like a $t$-th
exponential moment, and recall that, for any random variable $Y$, the
moment $\EE [\ee^{ t Y}]$ runs like $\ee^{ t\, {\rm esssup} Y}$.)

Another fundamental result is on the almost sure behaviour of the
total mass, which reads as follows in the above special case (see
Theorem~5.1 and Remark~3.17 in \cite{WK-K16}).

\begin{theorem}[Almost-sure asympotics of the PAM]
  Under the same assumptions as in Theorem~\ref{WK-ThmMomAsy}, there
  is a number $\widetilde\chi\in[0,2d]$ such that
\begin{equation}\label{WK-asTotMass}
  \frac{1}{t}\log U(t)=h^{}_t -\widetilde \chi+o(1),\qquad
  t\to\infty, \mbox{ almost surely},
\end{equation}
where $h_t=\varrho\log\log |B_t|$, and $|B_t|$ is the cardinality of
the box $B_t= [-t,t]^d\cap \ZZ^d$.
\end{theorem}

Indeed, $h_t$ is the asymptotics of the maximum of all the potential
values in $B_t$, and $\widetilde \chi$ is another variational formula
that describes the optimal principal eigenvalue of the operator
$\Delta+q$ under the assumption that the potential
$q\colon\ZZ^d\to \RR$ is not too improbable to be realised
under $\xi$. The explanation is that, in a box of some radius $r_t$
around the starting site $0$, one searches for a local region in which
the local principal eigenvalue of $\Delta+\xi$ is as large as
possible. Then the path in the Feynman--Kac formula runs quickly to
that place (which is typically $\approx r_t$ away) and stays in that island
for the remaining time until $t$. Here the radius $r_t$ is
chosen in such a way that an optimal compromise is reached between the
probabilistic costs to reach that remote place in a short time and to
have enough time left to benefit from that favourable potential.
As a deeper investigation in \cite{WK-BKS18} has shown, the
optimal choice is $r_t=\varrho t/\log t\log_3 t$ (with $\log_k t$ the
$k$-fold iteration of $\log$), and the time that is used for getting
there is of the order $t/\log t \log_2 t\log_3 t$.

Like for the moments, it is characteristic for the double-exponential
distribution in \eqref{WK-DE} that the size of the intermittent island
does not depend on $t$. Other potential distributions lead to
different sizes of the intermittent islands. As a rule, the heavier
the upper tails of the potential variable is, the smaller the island
is. E.g., for standard normal, the Weibull\index{distribution!Weibull}
or the Pareto\index{distribution!Pareto} distribution,
one obtains single-site islands, while for bounded random variables
the islands depend on $t$ and grow unboundedly.

\section{Multitype branching random walk in random
  potential}\label{WK-sec:MultiType}

In this section, we review the results of \cite{WK-GKS15}
about a version of the spatial branching walk
model\index{branching!process!multitype} defined in Section~\ref{WK-sec:model}
with the additional feature that we give \emph{two} categories of properties to
each particle, a spatial information and information about its
type. We now interpret $\cX$ as the location space and take
some Markov chain on it as a basis, and concerning the type, we
introduce an additional discrete set $\cT$, the \emph{type space}.
We will let the rates depend on the types of the particle
and of the offspring. Mathematically, there are certainly
possibilities to conceive $\cX \times \cT$ as one state
space, but we are interested in the effect coming from the difference
between the two components, since there are different modeling ideas
behind $\cX$ and $\cT$.

We consider a natural and flexible discrete-time version. Denote the
transition matrix of the Markov chain on $\cX$ by
$P=(P_{xy})_{x,y\in\cX}$. We equip $\cT$ with a set $\cA$ of directed
edges $(i,j)\in\cT\times\cT$ and obtain a directed finite graph
$\cG =(\cT ,\cA )$. We assume that each directed edge appears at most
once in $\cA$, and for each $i\in\cT$, there is at least one $j\in\cT$
such that $(i,j)\in\cA$. Self-edges $(i,i)$ may appear in $\cA$. To
each $y\in\cX$ we attach a matrix
$F_y=(F_y^{\ssup{i,j}})_{(i,j)\in\cA}$ of probability distributions on
$\NN_0$, the environment.  Given $F=(F_y)_{y\in\cX}$, we define a
discrete-time Markov process $(\eta_n)_{n\in\NN_0}$ on
$\NN_0^{\cT\times \cX}$, where $\eta_n(i,x)$ is the number of
particles of type $i$ at site $x$ at time $n$. The environment $F$
does not depend on time and is fixed throughout the evolution of
particles. We specify the transition mechanism of
$(\eta_n)_{n\in\NN_0}$ as follows: for a given state $\eta$, during
the time interval $(n,n+1)$, given that the configuration is equal to
$\eta$ at time $n$,
\begin{enumerate}
\item A particle of type $i$ located at site $y\in\cX$ produces,
  independently for $j\in\cT$ such that $(i,j)\in\cA$, precisely $k$
  particles of type $j$ at the same site $y$ with probability
  $F_y^{\ssup{i,j}}(k)$, for any $k\in\NN_0$. This offspring
  production is independent over all the particles in $\cX$ and over
  the time $n\in\NN_0$.

\item Immediately after creation, each new particle at $x$
  independently chooses a site $y$ with probability $P_{xy}$ and moves
  there.

\item The resulting particle configuration is $\eta_{n+1}$.
\end{enumerate}

Fix a site $y\in\cX$ and type $j\in\cT$ as starting sites. We start
the Markov chain $(\eta_n)_n$ with the initial configuration
$\eta_0(i,x)=\delta_{j}(i)\delta_{y}(x)$, and by ${\tt P}_{j,y}$ and
${\tt E}_{j,y}$ we denote its distribution and expectation,
respectively.  Note that the Markov chain depends on the realisation
of the environment $F$.  We are interested in the expectation of the
global population size,
$|\eta_n|:=\sum_{i\in\cT,x\in\cX}\eta_n(i,x)$,
\begin{equation}
  u^{}_n(i,x):={\tt E}^{}_{i,x}[|\eta^{}_n|],\qquad n\in\NN^{}_0, x\in\cX,
  i\in\cT .
\end{equation}

As the first main result, we develop two formulas for $u_n$: as the
$n$-th power of a certain $\cX\times\cT$-matrix, and as an expectation
of a multiplicative functional over $n$ steps of a particular Markov
chain on that set. These two formulas are extensions of well-known
representations of the expected number of particles of multitype
branching processes without space; we like the fact that one can see
these formulas from the viewpoint of a Feynman--Kac formula.

To formulate this, we need to introduce a Markov chain
$T=(T_n)_{n\in\NN_0}$ on the type space $\cT$ with transition
probabilities
\begin{equation}
p^{}_{ij}=\frac{\one\{(i,j)\in\cA\}}{{\rm deg}(i)},\qquad i,j\in\cT ,
\end{equation}
where ${\rm deg}(i)=|\{k\in \cT \colon (i,k)\in \cA \}|$ is the
  outdegree of $i$. We define $T$ and $X$ independently on a common
  probability space and write $\PP^{\ssup {T,X}}_{i,x}$ and
  $\EE^{\ssup{T,X}}_{i,x}$ for probability and expectation,
  respectively, where $T$ starts from $i$ and $X$ from $x$. We denote
  by $m_{ij}(y)=\sum_{k\in\NN_0}kF_y^{\ssup{i,j}}(k)$ the expectation
  of $F_y^{\ssup{i,j}}$ (the offspring expectation) and collect these
  random numbers in the matrix $M(y):=(m_{ij}(y))_{i,j\in\cT}$, where
  we put $m_{ij}(y)=0$ if $(i,j)\notin\cA$. The following is
  \cite[Proposition 1]{WK-GKS15}.

\begin{prop}[Representations for $u_n$]\label{WK-Prop-FKrepr}
For any $i\in\cT$ and any $x\in\cX$ and any $n\in\NN_0$, 
\begin{eqnarray}
  u^{}_n(i,x) &=&   \EE^{\ssup{T,X}}_{i,x}\Big[  \prod_{l=1}^{n}
                  \big( m^{}_{T_{l-1}T_l}(X^{}_{l-1})
                  {\rm deg}(T^{}_{l-1}) \big) \Big], \label{WK-AnneMom2} \\
  u^{}_n(i,x)&=&B^n\one (i,x)=\sum_{j\in\cT},
                 y\in\cX B^n_{(i,x),(j,y)},\label{WK-AnneMom2a} 
\end{eqnarray}
where $B$ is the $(\cT\times\cX )\times (\cT\times\cX )$ matrix with
coefficients
\begin{equation}\label{WK-Bdef}
B^{}_{(i,x),(j,y)}=m^{}_{ij}(x)P^{}_{xy}\one\{(i,j)\in \cA\} .
\end{equation}
\end{prop}

The equation in \eqref{WK-AnneMom2a} easily follows from
\eqref{WK-AnneMom2} by writing out explicitly the expectation over the
Markov chain and the $n$-fold matrix product. The interpretation of
\eqref{WK-AnneMom2} is as follows. Every $n$-step path
$((X_0,T_0),\dots,(X_n,T_n))$, together with the product over the
$m_{T_{l-1}T_l}(X_{l-1})$, stands for the union of all $n$-step
branching subtrees that produce, in the $l$-th step, for any
$l\in\{1,\dots,n\}$, from a particle of type $T_{l-1}$ located at
$X_{l-1}$ a particle of type $T_l$ that makes a step to $X_l$
immediately after creation. The product over the
$P_{X_{l-1}X_l}p_{T_{l-1}T_l}$, together will the product over the
$\deg(T_{l-1})$ (note the partial canceling of terms), summarises the
probabilities of all the jumps in the state space. In this way, we
encounter a discrete-time version of a Feynman--Kac formula for the
Markov chain $(T,X)$ on $\cT\times \cX$, however
with an interesting difference: the potential $\log m_{i,j}(x)$
depends on the \textit{vertices} of the space $\cX$ and on the
\textit{edges} of the type space $\cT$.

So far, the environment $F$ was fixed, but let us now turn to the case
of a \emph{random} environment\index{environment, random}, and let us describe our
assumptions. We assume that the collection of all distributions
$F_y^{\ssup{i,j}}$ with $y\in\cX$ and $(i,j)\in\cA$ is
independent. Their distribution depends on $(i,j)$, but not on $y$. We
call $F=(F_y)_{y\in\cX}$ the \textit{random environment} and denote by
${\rm Prob}$ and $\langle\cdot\rangle$ probability and expectation
with respect to $F$, respectively. Since we are only interested 
in the expectation of the global number of particles here, we will make our
assumptions on the environment only in terms of the quantities
$m_{ij}(y)=\sum_{k\in\NN_0}kF_y^{\ssup{i,j}}(k)$. In particular, we
assume that the collection of the $m_{ij}(y)$ is stochastically
independent over $y\in\cX$ and $i,j\in\cT$ (but of course not
identically distributed).

We will study the case where the upper tails of $ m_{ij}(y)$ lie in
the vicinity of the \textit{Weibull distribution},
\index{Weibull distribution}\index{distribution!Weibull} with parameter
$1/\rho_{ij}\in(0,\infty)$, i.e.,
    \begin{equation}
      {\rm Prob}(m^{}_{ij}(y)>r)\approx\exp\{ -r^{1/\rho_{ij}}\}, \qquad
      r\to\infty. \label{WK-DExpTail}
    \end{equation}
    Hence, $\log m_{ij}(y)$ lies in the vicinity of the
    double-exponential distribution defined in \eqref{WK-DE} with
    parameter $\rho_{ij}$. Let
    $H_{ij}(t):= \log\langle m_{ij}(y)^t\rangle$ denote the logarithm
    of the moment generating function. For
    $(i,j)\in\cT^2\setminus \cA$, we put
    $\rho_{ij}=0$. Hence, our environment distribution is
    characterised by the matrix-valued parameter
    $\rho=(\rho_{ij})_{i,j\in\cT}$. The larger $\rho_{ij}$
    is, the thicker are the tails of $ m_{ij}(y)$, i.e., the easier it
    is for $m_{ij}(y)$ to achieve extremely high values.

    The second main result of~\cite{WK-GKS15} is a formula for the
    large-$n$ asymptotics for the expectation of the Feynman--Kac
    formula from~\eqref{WK-AnneMom2}. This formula is
    in the spirit of~\eqref{WK-MomAsy}, however, it has some
    pecularities and some novelties. Let us first mention that there
    are basically two possible lines of proof for deriving the result:
    one about the large deviations for the empirical pair measure of
    the underlying Markov chain, and one using Frobenius eigenvalue
    theory. We decided to carry out only the former one. For any
    discrete set $S$, we denote by $\cM_1(S)$ the set of
    probability measures on $S$ and by
    $\cM_1^{\ssup {\rm s}}(S^2)$ the set of probability
    measures on $S^2$ with equal marginals. The first quantity of
    interest is
\begin{equation}\label{WK-lambdadef}
 \lambda(\rho) = \sup \big\{ \langle\mu,\rho\rangle\colon 
 \mu \in \cM_1^{\ssup {\rm s}}(\cT^2) \big\},\qquad
 \mbox{where }\langle\mu,\rho\rangle=\sum_{(i,j)\in\cA } \mu(i,j)\rho^{}_{ij},
\end{equation}
and the set of the corresponding maximisers:
\begin{equation}\label{WK-Lambdadef}
  \Lambda(\rho):=\big\{\mu\in\cM_1^{\ssup {\rm s}}(\cT^2)\colon
  \langle\mu,\rho\rangle=\lambda(\rho)\big\}.
\end{equation}
We introduce some notation. Each measure
$\nu\in\cM_1^{\ssup {\rm s}}((\cT\times\cX)^2)$ has a number of
marginal measures that are defined on different spaces, but in order
to keep the notation simple, we denote by $\overline{\nu}$ all these
marginals, namely,
\begin{equation}\label{WK-projections}
\begin{aligned}
  \overline{\nu}(i,j,x)&=\sum_{y\in\cX}\nu((i,x),(j,y)),\qquad
  \overline{\nu}(i,x)=\sum_{j\in\cT}\overline{\nu}(i,j,x),\\
  \overline{\nu}(i,j)&=\sum_{x\in\cX}\overline{\nu}(i,j,x),\qquad
  \overline{\nu}(i)=\sum_{j\in\cT}\overline{\nu}(i,j).
\end{aligned}
\end{equation}

To describe the second term in the asymptotics, we need to introduce
two functionals on probability measures
$\nu \in \cM_1^{\ssup s} (( \cT\times\cX )^2)$, an energy functional
$\cS$ and an entropy functional $\cI$. Indeed, define
\begin{eqnarray*}
  \cS (\nu)&:=&\sum_{(i,j)\in\cA} \rho^{}_{ij}\sum_{x\in\cX}
                \overline{\nu}(i,j,x)\log \overline{\nu}(i,j,x)+
   \sum_{(i,j)\in\cA}\overline{\nu}(i,j) \rho^{}_{ij}\log \rho^{}_{ij},\\
  \cI (\nu) &:=& \sum_{(i,j)\in\cA}\sum_{x,y \in \cX} \nu \big( (i,x),(j,y)\big)
                 \log\frac{\nu\big( (i,x),(j,y)\big)}{\overline{\nu}(i,x)
                 P^{}_{xy} }.
    \end{eqnarray*}
    We set $\cI (\nu)=\infty$ if $\nu $ is not absolutely continuous
    with respect to the measure
    $((i,x),(j,y))\mapsto \overline{\nu}(i,x)
    P_{xy}1\{(i,j)\in\cA\}$. Then $\cI (\nu)$ is equal to the relative
    entropy of $\nu$ with respect to this measure, up to the missing
    normalisation; note that the reference measure is not normalised,
    but has mass equal to $\sum_{i\in\cT}{\overline \nu}(i)\deg(i)$.

Now we can state our main result, \cite[Theorem 3]{WK-GKS15}: 

\begin{theorem}\label{WK-thm-Main} For any $i\in\cT$ and
  $x\in\cX$, as $n\to\infty$,
    \begin{equation}\label{WK-AnneMomAsy}
      \langle u^{}_n(i,x) \rangle = (n!)^{\lambda(\rho)}
      \ee^{-n\chi(\rho)} \ee^{o(n)}=  \exp
\Big(\lambda(\rho)n \log \frac{n}{\ee} - n \chi(\rho) + o(n)\Big),
    \end{equation}
where
\begin{equation}\label{WK-Chi}
\chi(\rho) = \inf \big\{ \cI (\nu)-
\cS (\nu) \colon \nu \in \cM_1^{\ssup s}\big( ( \cT \times \cX )^2 \big),
\overline{\nu}\in \Lambda(\rho)\big\}.
    \end{equation}
\end{theorem}

The central object in the proof and in the understanding of this
result is the \textit{empirical pair measure}
\begin{equation}\label{WK-LocTimes}
     \nu^{}_n = \frac{1}{n} \sum_{l=1}^n \delta^{}_{((T^{}_{l-1},X^{}_{l-1}),(T^{}_l,X^{}_l))}.
\end{equation}
In terms of the space-time random walk $(X,T)$, the number
$n \nu_n((i,x),(j,y))$ is equal to the number of $j$-type offspring of
any $i$-type particle located at $x$ by time $n$ that makes a step to
$y$ right after creation. Hence, $\nu_n$ stands for the union of all
$n$-step paths $((X_0,T_0),\dots,(X_n,T_n))$ that make precisely
$n \nu_n((i,x),(j,y))$ steps $(i,x)\to (j,y)$ for every
$i,j\in\cT$ and every $x,y\in\cX$. The term
$\cI (\nu)$ is the negative exponential rate of the probability
of this union under the Markov chain $X$, together with the
combinatorial complexity of the trajectories of types, and
$\cS (\nu)$, together with the leading term $\lambda(\rho)$,
is the one under the expectation w.r.t.~the random environment.

Theorem~\ref{WK-thm-Main} in particular shows that the main
contribution to the annealed moments of the numbers of particles,
$\lambda(\rho)$, comes from those $n$-step branching process subtrees
that produce, for some $\mu\in\Lambda(\rho)$, at approximately
$n\mu(i,j)$ of the $n$ steps a number of $j$-type particles from one
or more $i$-type particles, for any $i,j\in \cT$. Then the
value $\langle\mu,\rho\rangle$ gives the leading contribution on the
scale $n\log\frac{n}{\ee}$. It is interesting to note that the
optimality of the leading term has nothing to do with the spatial part
of the branching process, but only with the creation of particles. The
reason is that all the probabilities of spatial actions, i.e., of the
random walk $X$, are on the scale $n$, but the values of the offspring
expectations $m_{ij}(x)$ are typically on the scale
$\ee^{O(\log n)}$.

In this light, let us analyse the leading term $\lambda(\rho)$ a bit
more closely. A \textit{simple cycle} on $\cG$ is a path
$\gamma=(i_1,\dots,i_l)$ in $\cT$, with steps $(i_{m},i_{m+1})$ in
$\cA$, that begins and ends at the same vertex $i_1=i_{m+1}$, but
otherwise has no repeated vertices or edges. We write $(i,j)\in\gamma$
if the directed edge $(i,j)$ belongs to $\gamma$, that is, if
$(i,j)=(i_m,i_{m+1})$ for some $m\in\{1,\dots,l\}$. We call
$|\gamma|=l$ its \textit{length}. We denote by $\Gamma_l$ the set of
all simple cycles of length $l$ and by $\Gamma$ the set of all simple
cycles. We define the uniform measure on the edges of a simple cycle
$\gamma$,
\begin{equation}\label{WK-nugammadef}
\mu^{}_{\gamma}(i,j)=
\begin{cases}1/|\gamma|&\text{if }(i,j)\in \gamma,\\
0&\text{otherwise}.
\end{cases}
\end{equation}
It is clear that $\mu_\gamma\in\cM_1^{\ssup{\rm s}}(\cT^2)$ for any
$\gamma\in \Gamma$. Simple cycles are important for the asymptotics of
the annealed moments because the set of extreme points of
$\cM_1^{\ssup {\rm s}}(\cT^2)$ consists exactly of the $\mu_\gamma$'s
with $\gamma$ the simple cycles of the graph $\cG$. The following is
\cite[Lemma 1]{WK-GKS15}.

\begin{lemma}\label{WK-extrarecycles}
  The set of extreme points of the convex set
  $\cM_1^{\ssup {\rm s}}(\cT^2)$ is equal to
  $\{\mu_\gamma\colon\gamma\in\Gamma\}$.
\end{lemma}

Since the optimisation problem in (\ref{WK-lambdadef}) is a linear
optimisation problem on the convex, compact set
$\cM_1^{\ssup{\rm s}}(\cT^2)$, the Krein-Milman theorem and Lemma
\ref{WK-extrarecycles} imply the following characterisation of the
leading term in \eqref{WK-AnneMomAsy}. The following is \cite[Lemma
2]{WK-GKS15}.

\begin{lemma}\label{WK-lem-MaxCyclesLemma}
$$
\lambda(\rho) =  \max\Big\{\langle \nu^P_{\gamma},\rho\rangle
  :\gamma\in\Gamma\Big\} =\max\Big\{ \frac{1}{|\gamma|} \sum_{m=1}^{|\gamma|}
  \rho^{}_{i^{}_{m-1}i^{}_m}\colon (i^{}_1,\dots,i^{}_{|\gamma|})\in \Gamma \Big\}.
$$
\end{lemma}

The interpretation of Lemma~\ref{WK-lem-MaxCyclesLemma} is that the
leading contribution to the expected population size comes from
optimal cycles $(i_1,\dots,i_{|\gamma|})\in\Gamma_{|\gamma|}$. In
terms of branching process trees, they are considered optimal if they
produce only $i_{m+1}$-type particles from $i_m$-type particles for
any $m\in\{1,\dots,|\gamma|\}$ (with $i_{|\gamma|+1}=i_1$), but no
other offspring.

\section{The PAM on finite graphs}
\label{WK-sec:hypercube}

Another, more biologically inspired, direction in the
investigation of the PAM, is the PAM on some finite graph (like the
full graph or the hypercube) and the concentration on the question
about the amount of time that one has to give to the branching process
such that the overwhelming part of the population has found the way to
the ``fittest'' site. Here we rely on the interpretation that we
explained in Section \ref{WK-sec:MutSel}, and the question is
investigated in the limit of a large (but finite) graph and late
times, and the relation of the growths of time and space is
crucial. In order to make a mathematical treatment feasable without
too many technicalities, one typically assumes the random potential
(the ``fitness landscape'')\index{fitness!landscape} to be i.i.d.\
exponentially distributed times a growth parameter, such that one is
in the regime where the intermittent islands (in the understanding of
Section~\ref{WK-sec:MomAsy}) are singletons.

The first work on the PAM on a finite graph
$\cX$\index{parabolic Anderson model!on hypercube} was --- to
the best of our knowledge --- \cite{WK-FM90}, which considered the
complete graph with $N$ nodes and an exponentially distributed i.i.d.\
random potential $\xi$. However, instead of $\Delta$ in
\eqref{WK-PAM}, the \emph{rescaled\/} Laplace operator $\frac{1}{N}\Delta$
is considered. This is mathematically equivalent to \eqref{WK-PAM}
with graph-size dependent potential $N\xi$; the large prefactor $N$
supports a strong concentration of the mass of branching particles in
small intermittent islands; actually here we are concerned with single
sites.

In that work, the initial condition was taken to be localised in the
site of the $k$-th largest of the potential values, and the question
was raised, for what choices of $t=t_N$, in the large-$N$ limit, the
main mass of the particle system travels to the site of the largest
potential value, and for what choices it stays at the initial site by
time $t$. The authors found the leading scale and the criterion for the
answer, and they derived the first term in the asymptotics for
the expectation of the total mass with this initial condition. We do 
not go deeper into these details.

In \cite{WK-AGH16}, the same question was raised for the hypercupe
$\{-1,1\}^N$ with an i.i.d.\ potential, the assumptions of which are
formulated by requiring that $2^N$ i.i.d.\ copies of the potential
variables leave gaps of order one between their consecutive leading
order statistics, asymptotically as $N\to\infty$. According to
standard extreme-value statistics, this includes the case of centred
Gaussians with variance $N$. The main result is the following (see
\cite[Theorem 1.2]{WK-AGH16}).

\begin{theorem}\label{WK-ThmHypcubeAsy}
  Assume that $U_N(t)$ is the total mass of the solution to the PAM in
  \eqref{WK-PAM} on $\cX =\{-1,1\}^N$ with $\Delta$ replaced by
  $\frac{1}{N} \Delta$ and with an i.i.d.\ random potential $\xi$ on
  $\cX$ such that their leading order statistics
  $\max_\cX\xi= \xi_{1,2^N}>\xi_{2,2^N}>\xi_{3,2^N}>\dots$
  leaves gaps of order one between its consecutive values,
  asymptotically as $N\to\infty$. We assume that the initial condition
  in \eqref{WK-PAM} is the delta-measure in the site of $\cX$
  in which $\xi_{k,2^N}$ sits with a fixed integer $k>1$. Fix
  $\varepsilon>0$. If
  $t_N\geq (1+\varepsilon)(N\log N)/2(\xi_{1,2^N}-\xi_{k,2^N})$, then,
  almost surely,
$$
U^{}_N(t)=\exp\Big\{(\xi^{}_{1,2^N}-1)t_N+
\frac{1}{2} (1+o(1))N\log N\Big\},\qquad N\to\infty.
$$
In contrast, if
$t_N\leq (1-\varepsilon)(N\log N)/2(\xi_{1,2^N}-\xi_{k,2^N})$, then,
almost surely,
$$
U^{}_N(t)=\exp\{(\xi^{}_{k,2^N}-1)t^{}_N\}\big( 1+o(1)\big) ,\qquad N\to\infty.
$$
\end{theorem}

The interpretation is that the main mass has, in the first case of a
large time horizon, enough time to move to the site of the maximal
potential value, $\xi_{1,2^N}$, while in the second case (small time
horizon), it stays in the initial site.

\section{Higher moments of the numbers of
  particles}\label{WK-sec:HigMom}

As we mentioned in Section \ref{WK-sec:PAM}, the first
moment (the expectation) $u(t,x)$ of the number of particles in the branching
random walk system at time $t$ in the site $x$, satisfies
the heat equation with random potential in \eqref{WK-PAM}, which is
rather amenable to a deeper analysis, due to the Feynman--Kac formula
in \eqref{WK-FK} and the eigenvalue expansion in
\eqref{WK-eigenexp}. However, this concerns only the \emph{first}
moment, and this is quite short of a deeper understanding of the
entire particle system. Substantially more information is contained in
the $n$-th moments of the particle number at time $t$ at site $x$ for
$n\in\NN$, i.e., in the expectation of $\eta(t,x)^n$, again
taken only over the migration, branching and killing mechanisms with
fixed branching/killing rates. In this section, we are going to report
on the work \cite{WK-GKS13} on this aspect. To keep things simple in
this survey, we only consider here $U_n(t)$, the expectation of
$(\sum_{x\in\cX}\eta(t,x))^n$. Hence, the case $n=1$ is the one
that we handled before. Again, like in Section \ref{WK-sec:MultiType},
we will present two main results: An explicit formula of Feynman--Kac
type for $U_n(t)$ for fixed branching/killing rates, and the
identification of the large-$t$ asymptotics of the expectation of
$U_n(t)$, taken with respect to a particular choice of the
distribution of these random rates.

One of our motivations to study $U_n(t)$ is that asymptotic knowledge
on the behaviour of $\langle U_n(t)\rangle $ might be able to tell
something about the distributional behaviour of the total number of
particles, as in the general theorem where (under some technical
conditions) the convergence of all the $n$-th moments of a random
variable implies its weak convergence towards a random variable whose
$n$-th moments are equal to the limits.\index{branching!process!higher moment}

Meanwhile, also deep investigations have been carried out directly on
the branching particle system in the spirit of the PAM and using a
great deal of the results and methods specific for the PAM \cite{WK-OR18},
but this concerned only the simplest potential distribution (the
Pareto distribution) and was technically enormously cumbersome.

The method of choice to analyse $U_n(t)$ is to derive a similar
equation as the heat equation with random potential and to try to
employ the techniques that were helpful in the analysis of the PAM. In
earlier work \cite{WK-ABMY00}, a recursive equation was derived for
$u_n(t,x)$ (defined as the expectation of $\eta(t,x)^n$) in terms of
all the functions $u_1,\dots,u_{n-1}$. This formula was explicit
enough to derive, for a particular potential distribution that leads
to single-site islands, the first term in the large-$t$ asymptotics of
the moments of $u_n(t,\cdot)$. However, this is obviously
unsatisfactory for several reasons. Instead, in \cite{WK-GKS13}, a
direct formula for the expectation of $u_n$ (in fact, for its $p$-th
moment for any $p\in\NN$) is derived, and this is so explicit
that, for the double-exponential distribution in
\eqref{WK-DE}\index{distribution!double-exponential}, the two main
terms could be derived in terms of a variational formula that admits
an interpretation and further analysis. The main tool for deriving the
direct formula is the many-to-one formula \index{many-to-one formula}from the theory of branching processes, derived via spine
techniques.\index{spine technique}

Unlike in the case $n=1$ of the first moment, the $n$-th moment is not
just a function of the difference $\xi=\xi_+-\xi_-$, but we need to
keep track of both the splitting rate $\xi_+=\xi_2$ and the killing
rate $\xi_-=\xi_0$.

Let us state the formula from \cite[Theorem 2.1]{WK-GKS13} in the
special case $n=2$.

\begin{theorem}[Feynman--Kac-like formula for $U_2$]
  Let $U_n(t)$ denote the $n$-th moment of the total number of
  particles, for any $n\in\NN$. Then, in the case $n=2$, we have
  $U_2=U_1+\widetilde U_2$, where
\begin{equation}\label{WK-secondmom}
\widetilde U^{}_2(t)=\int_0^t \EE^{}_\cO\Big[{\rm e}^{\int_0^s \xi(X^{}_r)\,{\rm d} r}{\rm e}^{\int_s^t \xi(X_r')\,{\rm d}r}{\rm e}^{\int_s^t \xi(X_r'')\,{\rm d}r}2\xi^{}_2(X^{}_s)\Big]\,{\rm d} s,
\end{equation}
where $(X_r)_{r\in[0,s]}$, $(X'_r)_{r\in[s,t]}$ and
$(X''_r)_{r\in[s,t]}$ are three independent simple random walks,
conditioned on $X_0=\cO$ and $X'_s=X''_s=X_s$.
\end{theorem}

In other words, the right-hand side of \eqref{WK-secondmom} is the
expectation over a branching random walk with precisely one splitting
event in the time interval $[0,t]$, namely at time $s$. The first term
in the representation $U_2=U_1+\widetilde U_2$ corresponds to absence
of splitting, the second to exactly one splitting. For general $n$,
the formula for $U_n$ is similar (one has to sum over all numbers of
splitting events in $\{0,\dots,n-1\}$), but much more cumbersome,
since combinatorial prefactors are involved and multiple powers of
$\xi_2$-values at the splitting locations. The formula in
(\ref{WK-secondmom}) is the main result of \cite{WK-GKS13} about a
representation of the $n$-th moments of the number of particles, for
fixed branching/killing rates, in the special case $n=2$.

Now we turn to the second question. For the analysis of the large-$t$
asymptotics of the expectation of $\widetilde U_2(t)$ and for its
intuitive understanding, the formula in \eqref{WK-secondmom} is well
suitable, at least for the case that both random fields $\xi_2$ and
$\xi_0$ are i.i.d.\ sequences of double-exponentially distributed
variables as in \eqref{WK-DE} and are independent. The first
observation is that the term $2\xi_2(X_s)$ should hardly have any
influence, and this is indeed true under some technical assumption
that forbids too thick upper tails of the random variable $\xi_2(0)$.

The second observation is that, as we are considering the case that
the potential can assume very large (positive) values, it should
be giving particularly large values if the splitting time $s$ is as
early as possible, since it is then \emph{two} independent random walks
that can contribute to the random-potential expectation: it is indeed
the expectation of a \emph{square} of the Feynman--Kac formula, since
each of the two copies contributes independently like one total
mass. This reasoning can be done already when considering the leading
terms of the expectation: Since we have three random walks with
running time $s$, $t-s$ and $t-s$, one should expect that the leading
logarithmic term is roughly equal to $H(s+(t-s)+(t-s))=H(2t-s)$, which
is maximal if $s$ is minimal, since we assume that
$\frac 1t H(t)\to\infty$ (as $\xi$ is unbounded). Hence, the first
term in the asymptotics for
$\frac 1t\log \langle \widetilde U_2(t)\rangle$ should be $H(2t)$, and
the second one should be $-2 \chi(p)$ in the notation of
\eqref{WK-MomAsy}. The same picture is true for the higher moments,
i.e., for the moment asymptotics of $U_n(t)$, with 2 replaced by $n$
at both occurrences. This has indeed been proved as the second main
result of \cite{WK-GKS13} in Theorem 1.3 there:

\begin{theorem}[Moment asymptotics for the branching random walk in
  random environment]
  Let the random potentials $\xi_2=(\xi_2(x))_{x\in\ZZ^d}$ and
  $\xi_0=(\xi_0(x))_{x\in\ZZ^d}$ be independent and i.i.d.\
  and double-exponentially distributed as in \eqref{WK-DE} with
  parameter $\rho\in(0,\infty)$. Denote by $U_n(t)$ the $n$-th moment
  of the total number of particles in the branching random walk in the
  random environment $(\xi_0,\xi_2)$. Then, for any $n,p\in\NN$,
  $$
\langle U^{}_n(t)^p\rangle=\exp\big\{H(npt)-npt\chi+o(t)\big\},
\qquad t\to\infty,
$$
where $\chi\in[0,2d]$ is the number appearing in Theorem \ref{WK-ThmMomAsy}.
\end{theorem}

The most important conclusion is that, at least for the potential
double-exponentially distributed, the main contribution to the
expected $p$-th power of the total mass of the $n$-th moment comes
from early splittings in the moment formula. This implies that
$$
\langle U^{}_n(t)^p\rangle =\langle U(t)^{np}\rangle {\rm e}^{o(t)}
=\langle U(tnp)\rangle {\rm e}^{o(t)},\qquad t\to\infty.
$$
Recall that we believe that this phenomenon comes from the
unboundedness of the branching rates $\xi(0)$, more precisely from the
super-linear behaviour of the leading term, the logarithmic moment
generating function $H(t)$. If the potential random variable $\xi(0)$
is not positive, but attains only strictly negative values, then we
expect that the opposite behaviour is crucial, i.e., a very late
splitting, and the result should be that
$\langle U_n(t)\rangle =\langle U_1(t)\rangle {\rm e}^{o(t)}$ as
$t\to\infty$. If the essential supremum of the potential is zero, then
we expect that deeper investigations are necessary and that, in some
cases, richer pictures may appear.

\section{Further perspectives}
\label{WK-sec:future}

\subsection{High-moment asymptotics}

In ongoing work, we are currently deriving the large-$n$
asymptotics of the $n$-th moment (taken over all randomnesses:
migration, branching, killing, and all the rates) of the number of
particles in the branching random walk model of Section
\ref{WK-sec:model} on a finite time interval in $\ZZ^d$, i.e.,
the asymptotics of $\langle U_n(1)\rangle$. This will say something
about the most probable way of the branching particle system to
produce as much offspring as possible over a finite time interval,
i.e., about the questions how high the potential values should be, how
often and how quickly after each other the splitting of the particles
occurs, and how many steps all the many paths make.  For deriving
this, we are exploiting the moment formula that we explained in
Section \ref{WK-sec:HigMom}.

\subsection{Branching random walks on (random) graphs}\label{WK-sec:BRWgraphs}

As we explained in Section~\ref{WK-sec:hypercube}, an
investigation of the PAM on (possibly finite) graphs is biologically
sound, but has not been done yet on a broad front; we are actually 
aware only of the two works \cite{WK-FM90} and \cite{WK-AGH16}. In
another ongoing investigation, the state space $\cX$ is taken as a
\emph{random} graph that is locally tree-like\index{parabolic Anderson model!on random graph}. The main examples are a
Galton--Watson tree\index{parabolic Anderson model!on Galton--Watson tree} with bounded degrees and the configuration model. The main
goal is --- for the random potential double-exponentially distributed
as in \eqref{WK-DE} --- to find the large-$t$ asymptotics for the total
mass $U(t)$ of the solution with high probability and to describe the
structure of the parts of the random graphs that give the main
contribution, i.e., of the intermittent islands. The two main
interests here are to understand in general the influence (1) of the
exponential structure of the tree, i.e., the fact that the volume of a
ball with diameter $r$ runs exponentially fast in $r$, and (2) of the
randomness of the graph structure.

In the long-term, we hope to answer also the questions that are
answered in the case of the state space $\cX=\{-1,1\}^N$ in
Theorem~\ref{WK-ThmHypcubeAsy}, but the infiniteness of the state
space, its randomness and the fact that the intermittent islands
are not single sites, but have some structure make the question about
the location and the shape of the intermittent islands a big
enterprise.

\subsection{Self-repellent random walk in random
  potential}\label{WK-sec:SAW}

As described in Section~\ref{WK-sec:MomAsy}, the
Feynman--Kac formula that represents the total mass of the solution to
the PAM actually displays a random walk in random potential, and the
large-time asymptotics are carried by those random walk paths that
find their way to an optimal local region in the potential. The random
walk in random potential feels an attractive force towards the
extremal regions of the potential.

In another ongoing work, we are investigating the effect of
an additional counter force: some additional self-interaction
that suppresses self-intersections of the random walk until time
$t$. In other words, we replace the free walk by the well-known
\textit{weakly self-avoiding}\index{random!walk!self-avoiding}
or \textit{self-repellent walk}, which
is given by an exponential weight of the form
$$
\frac{1}{Z^{}_{T,\beta}}\exp\Big\{-\beta \int_0^T\int_0^T \one\{X^{}_s=X^{}_t\}\,
\dd s \dd t\Big\},
$$
where $\beta\in(0,\infty)$ is a parameter, and $Z_{T,\beta}$ is the
normalising constant. In this model, the path in the Feynman--Kac
formula cannot spend too much time anymore in single sites, the
intermittent islands. The goal is to describe what it does instead.
If the underlying state space is the lattice $\ZZ^d$, then the
most obvious conjecture is that it will visit not only one of the
intermittent islands, but several ones after each other, even though
they are far away from each other. This is a random strategy,
depending via an opimisation problem on the limiting spatial
extreme-value order statistics of the potential. We think that this
strategy is indeed optimal for most values of the thickness parameter
of the potential distribution, but for some values this does not seem
to be true, and we currently have no clue how to decribe the typical
behaviour of the path properly.




\end{document}